\documentclass[12pt]{article}
\usepackage{amsmath,amssymb,amsthm, amscd}

\begin{document}

\begin{center}
{\bf A characterization of nilpotent orbit closures among symplectic singularities II}
\end{center}
\vspace{0.4cm}

\begin{center}
{\bf Yoshinori Namikawa}
\end{center}
\vspace{0.2cm}

This short note is a supplement to the article [Na-Part I] with the same title.  
In [Na-Part I] we characterize 
a normal nilpotent orbit closure of a complex semisimple Lie algebra as a conical symplectic 
variety of maximal weight $1$. In this article we treat a conical symplectic variety obtained as a finite covering of a (not necessarily normal) nilpotent orbit closure of a complex semisimple 
Lie algebra. The normalization of a non-normal nilpotent orbit is among them. Before stating the result we prepare some notation. Let $R$ be a normal domain of finite type over  $\mathbf{C}$, where $R$ is positively graded $\mathbf{C}$-algebra $\oplus_{i \geq 0}R_i$ with  
$R_0 = \mathbf{C}$. Then $X = \mathrm{Spec}R$ is a conical symplectic variety if $X$ admits a holomorphic symplectic 2-form $\omega$ on the regular locus $X_{reg}$ such that $\omega$ 
is a homogeneous 2-form for the $\mathbf{C}^*$-action and extends to a regular 2-form on a resolution $\tilde{X} \to X$. We denote by $wt(\omega)$ the degree of the 
homogeneous 2-form $\omega$. We already know that $wt(\omega) > 0$ (cf.   
[Na 2], Lemma (2.2)). For a positive integer $l$ we denote by $R(l)$ the subring of $R$ defined as the image of the map $\oplus \mathrm{Sym}^i(R_l) \to R$. If we can write $R = \oplus_{i \geq 0} R_{ik}$ for some $k$, then we put $R_{[j]} := R_{jk}$ and introduce a new grading $R = \oplus_{j \geq 0}R_{[j]}$. Such a grading is called a {\em reduced grading} of $R$. 
Our main result is summarized as \vspace{0.2cm}

{\bf Theorem}. {\em Let $R$ be the coordinate ring of a conical symplectic variety $(X, \omega)$ 
with $wt(\omega) = l$. Assume that $R$ is a finite $R(l)$-module. Then} \vspace{0.15cm} 

(i) {\em $R_l$ is a semisimple Lie algebra $\mathfrak{g}$ and the simply connected algebraic group $G$ with $\mathrm{Lie}(G) = \mathfrak{g}$ acts on $X$ as $\mathbf{C}^*$-equivariant symplectic automorphism.}  

(ii) {\em The inclusion $R(l) \to R$ induces a finite $G$-morphism $f: X \to \bar{O}$ onto a nilpotent orbit closure $\bar{O}$ of $\mathfrak{g}$. Moreover, $f$ is \'{e}tale in codimension one and $f^{-1}(O) \to O$ is an \'{e}tale $G$-covering. The symplectic form $\omega$ is obtained as 
the pull-back of the Kirillov-Kostant form $\omega_{KK}$ on $O$ by $f$.} \vspace{0.15cm} 

{\em When $l $ is an even number, one has $R = \oplus_{j \geq 0}R_{jl/2}$.  When 
$l $ is an odd number, one has $R = \oplus_{j \geq 0}R_{jl}$. In particular, if $l$ is even, then the grading of $R$ can be reduced so that $wt(\omega) = 2$ and $R$ is a finite $R(2)$-module (with respect to the reduced grading).  If  $l$ is odd, then the grading of $R$ can be reduced so that $wt(\omega) = 1$ and $R$ is a finite $R(1)$-module (with respect to the 
reduced grading). } 
\vspace{0.2cm} 

We will prove the first part of Theorem in Proposition 2 and will prove the latter part in Proposition 3. The essential part is the semi-simplicity of $\mathfrak{g}$ and it is based 
on Proposition 4 of [Na-Part I]. As the following example suggests, the condition that 
$R$ is a finite $R(l)$-module is indispensable. \vspace{0.15cm}

{\bf Example}. For an integer $n \geq 2$ we put $R := \mathbf{C}[x,y,z]/(xy - z^{n+1})$ and 
introduce a grading of $R$ by $(wt(x), wt(y), wt(z)) = (1, n, 1)$. Then $X := \mathrm{Spec} R$ 
admits a symplectic 2-form $\omega := Res(dx \wedge dy \wedge dz/xy - z^{n+1})$ and 
$(X, \omega)$ is a conical symplectic variety of dim $2$ with $wt(\omega) = 1$. In this case 
$R(1)$ and $R$ have the same quotient field, but $R$ is not a finite $R(1)$-module. 
The algebraic group 
$$ G = \{ \left( \begin{array}{cc}
a & 0 \\ 
b & 1   
\end{array}\right) \vert \; \;  a \in \mathbf{C}^*, \; b \in \mathbf{C} \}$$ acts 
on $X$ as $\mathbf{C}^*$-equivariant symplectic automorphisms by 
$$x \to ax, \; \; z \to bx + z, \; \; y \to a^{-1}\{y + \sum_{i = 0}^{n} \binom {n+1}{i+1} b^{i+1}x^i z^{n-i}\}.$$ The degree 1 part $R_1$ is a Lie algebra by the Poisson bracket on $R$, and it coincides with the Lie algebra $\mathfrak{g}$ of $G$, which is not even reductive.        
\vspace{0.15cm}   

We first treat the case $l = 1$ with the additional condition that $R$ and $R(1)$ have the same quotient field.   

{\bf Proposition 1}. {\em Let $R$ be the coordinate ring of a conical symplectic variety 
$(X, \omega)$. 
We denote by $R(1)$ the image of the natural homomorphism $\phi: \oplus Sym^i(R_1)\to R$. Assume that $R$ is a finite $R(1)$-module and two rings $R$ and $R(1)$ have the same quotient field. Then one has $wt(\omega) = 2$ or $wt(\omega) = 1$.     
In the first case, $(X, \omega)$ is isomorphic to (i) and in the second case $(X, \omega)$ is 
isomorphic to (ii)}:  

(i) {\em  $(\mathbf{C}^{2d}, \omega_{st})$ with $\omega_{st} = \Sigma_{1 \le i \le d} dz_i \wedge dz_{i + d}$},

(ii) {\em the normalization of a nilpotent orbit closure 
$(\bar{O}, \omega_{KK})$ of a complex semisimple Lie algebra $\mathfrak{g}$}.   
 
{\em Proof}. Put $W := \mathrm{Spec} R(1)$. Then the ring homomorphism $\oplus \mathrm{Sym}^i(R_1) \to R(1) \subset R$ induces a morphism $f: X \to W \subset 
\mathbf{C}^n$, where $n = \dim R_1$. The morphism $X \to W$ is the normalization map.  
  
We already know that $wt(\omega) > 0$. First we prove $wt(\omega) \leq 2$. 
In fact, if $wt(\omega) > 2$, then the Poisson bracket has degree $< -2$, which means 
that $\{R_1, R_1\} = 0$; hence, $\{R(1), R(1)\} = 0$. Therefore $W$ has a trivial Poisson 
structure. The morphism $X \to W$ is a Poisson morphism which is generically an isomorphism. 
This implies that the Poisson structure of $X$ is trivial, which is a contradiction. Now assume that $l := wt(\omega) = 2$; then we have a skew-symmetric form $\{ \; , \;\} : R_1 \times R_1 
\to R_0 = \mathbf{C}$. The skew-symmetric form determines a Poisson structure on 
$\mathbf{C}^n := \mathrm{Spec} (\oplus \mathrm{Sym}^i(R_1))$ and $W$ is a Poisson closed subscheme of $\mathbf{C}^n$.  

If the skew-symmetric form is degenerate, then we have a nonzero element $x_1 \in R_1$ 
such that $\{x_1, R_1 \} = 0$. Notice that $\phi (x_1) \ne 0$ because $\phi\vert_{R_1}$ is an injection. Then we have $\{\phi (x_1), R(1)\} = 0$. We prove that $\{\phi (x_1), R\} = 0$. 
For $y \in R$, we take a monic polynomial $h(t) = t^n + a_1t^{n-1} + ... + a_n$ with coefficients 
in $R(1)$ so that $h(y) = 0$ and $\mathrm{deg}(h)$ is minimal among those. Since $\{\phi(x_1), a_i\} = 0$ for all $i$, we can write   $$0 = \{\phi (x_1), h(y)\} = (ny^{n-1} + (n-1)a_1y^{n-2} + ... + a_{n-1})\cdot \{\phi (x_1), y\}. $$  By the minimality of $\mathrm{deg}(h)$, we see that $ny^{n-1} + (n-1)a_1y^{n-2} + ... + a_{n-1} \ne 0$. As $R$ is an integral domain, we have $\{\phi (x_1), y\} = 0$.  This implies that the effective divisor $D$ of $X$ 
defined by $\phi(x_1) = 0$ is a closed Poisson subscheme. In this case $D \cap X_{reg} 
\ne \emptyset$ and $D \cap X_{reg}$ is a closed Poisson subscheme of $X_{reg}$. But this  contradicts the fact that the smooth symplectic variety $X_{reg}$ does not 
have any closed Poisson subscheme except $X_{reg}$ itself.   
   
Therefore the skew-symmetric form must be non-degenerate. In this case $n$ is even, say 
$n = 2d$, and $W$ is a closed Poisson structure of the affine space $(\mathbf{C}^{2d}, \omega_{st})$ with the standard 
symplectic form. Since $(\mathbf{C}^{2d}, \omega_{st})$ has no closed Poisson subscheme 
except $\mathbf{C}^{2d}$ itself,  $W = \mathbf{C}^{2d}$. This is the case (i).

We next assume that $l = 1$. Then $R_1$ has a structure 
of a Lie algebra $\mathfrak{g}$. We regard its dual space $\mathfrak{g}^*$ as an affine 
variety, which is nothing but the affine space $\mathbf{C}^n$. Let $G$ be the adjoint group of $\mathfrak{g}$ (cf. [Pro], p. 86) and let $\mathrm{Aut}^{\mathbf{C}^*}(X, \omega)$ be the $\mathbf{C}^*$-equivariant automorphism group of $(X, \omega)$. An element of  $\mathrm{Aut}^{\mathbf{C}^*}(X, \omega)$ determines a linear automorphism of $R_1 = 
\mathfrak{g}$ and conversely this linear automorphism uniquely determines an automorphism of $W$. Since $X$ is the normalization of $W$, the automorphism of $W$ uniquely extends to an automorphism of $X$. In this way one can embed $\mathrm{Aut}^{\mathbf{C}^*}(X, \omega)$ in $GL(\mathfrak{g})$. Then Proposition 2 of [Na-Part I] holds as it stands; namely, $G$ is the identity component of $\mathrm{Aut}^{\mathbf{C}^*}(X, \omega)$ and $\mathfrak{g}$ has no center.   

In the next step we will need a slight modification of Kaledin's result [Ka, Theorem 2.5], which can be stated as 

``{\em Let $Z$ be a Poisson variety and assume that its normalization $Z^{norm}$ is a symplectic variety with respect to the natural Poisson structure induced from $Z$ (cf. [Ka 2]). Then $Z$ is holonomic }".   

Actually, the statement of [Ka Theorem 2.5] is claimed only when $Z$ itself is a symplectic variety. However, the proof can be applied even in this case. For the convenience of readers 
we will give a proof.  Let $Y \subset Z$ be an integral Poisson subscheme. 
We have to prove that the Poisson structure on $Y$ is generically non-degenerate. Let 
$\pi: \tilde{Z} \to Z$ be the canonical resolution in the sense that a) $\pi$ is an isomorphism 
over the smooth locus $Z_{reg}$, and b) every vector field $\zeta$ on $Z$ lifts to a vector field 
$\tilde{\zeta}$ on $\tilde{Z}$.  Let $\nu: Z^{norm} \to Z$ be the normalization map. Then 
$\pi$ factorizes as $\tilde{Z} \stackrel{\pi'}\to Z^{norm} \stackrel{\nu}\to Z$.  
Take a non-empty Zariski open subset $U$ of $Y$ so that $U$ is smooth and the map $\pi^{-1}(U) \to U$ is generically smooth on $\pi^{-1}(u)$ for every 
closed point $u \in U$. The Poisson vector $\Theta$ for $Z$ induces a skew-symmetric 
form $T_y^*Y \times T_y^*Y \to \mathbf{C}$ for a point $y \in U$. Assume that it 
is a degenerate form. Then there is a function $f$ on $Z$ such that $df \ne 0$ at $y$, but 
the Hamiltonian vector field $H_f := \Theta (df, \cdot)$ vanishes at $y$. Generically on 
$Z$, there is a symplectic 2-form $\Omega$ such that $df = i_{H_f}\Omega$. 
Since $\pi$ is the canonical resolution, $H_f$ lifts to a vector field $\tilde{H}_f$ on $\tilde{Z}$. 
Generically we have $\pi^*(df) = i_{\tilde{H}_f}\Omega$; but, since $\Omega$ is defined everywhere on $\tilde{Z}$, this equality has a sense everywhere on $\tilde{Z}$.  
We put $U^{norm} := \nu^{-1}(U)$. Then $\pi$ factorizes as $\pi^{1}(U) \stackrel{\pi'}\to 
U^{norm} \stackrel{\nu}\to U$. By Lemma 2.9 of [Ka], we can take an open set $U'$ of 
$U^{norm}$ so that on each connected component $V$ of the smooth part of $(\pi')^{-1}(U')$ we have $$\Omega\vert_V = (\pi')^*\Omega_0$$ for some 2-form $\Omega_0$ on $U'$. Notice 
that $\nu$ is a closed map and $\nu(U^{norm} - U')$ is a proper closed subset of $U$. Replacing $U$ by $U - \nu(U^{norm} - U')$, we may assume from the first that $U' = U^{norm}$. Moreover we shrink $U$ so that $U^{norm} \to U$ is an \'{e}tale map. In general $\Omega_0$ is not the pull-back of a 2-form on $U$. But, since $U^{norm}$ is \'{e}tale over $U$, 
one can take an {\em analytic} open neighborhood $U(y)$ of $y \in U$ in such a way that 
$\Omega_0\vert_{\nu^{-1}(U(y))} = \nu^*(\Omega_0)'$ for some complex analytic 2-form  
$(\Omega_0)'$ on $U(y)$. Now we have an equality on $\pi^{1}(U(y))$: 
$$i_{\tilde{H}_f}\Omega = \pi^*( i_{H_f}(\Omega_0)').$$ 
Since $H_f$ vanishes at $y$, the right hand side vanishes on $\pi^{-1}(y)$. On the other 
hand, left hand side is $\pi^*(df)$, which does not vanish on $\pi^{-1}(y)$ because $\pi^{-1}(U) 
\to U$ is a generically smooth on $\pi^{-1}(y)$. This is a contradiction; hence the Poisson scheme $Y$ is generically non-degenerate. 
\vspace{0.2cm}
     
Now that $Z$ is holonomic,   
we see that there are only finitely many integral Poisson closed  subschemes in $Z$ by [Ka, Proposition 3.1]. Let us return to the proof of Proposition 1. Since $G$ is the identity component of 
$\mathrm{Aut}^{\mathbf{C}^*}(X, \omega)$, the coadjoint action of $G$ on $\mathfrak{g}^*$ 
preserves $W$ as a subvariety. Applying the result above to $W$, we can prove in a similar way to [Na-Part I, Proposition 3] that $W$ is the closure of a coadjoint orbit $(O, \omega_{KK})$ of $\mathfrak{g}^*$. 
Finally we must show that $\mathfrak{g}$ is semisimple. By the observation just above, 
$W$ has only finitely many integral Poisson closed subschemes. On the other hand, if $\mathfrak{g}$ is not semisimple, then $(\bar{O}, \omega_{KK})$ must 
have infinitely many such subschemes by [Na-Part I, Proposition 4].  This is a contradiction. 
Q.E.D.  \vspace{0.2cm}

{\bf Remark}. Let $V := \mathrm{Spec}\; S$ be an affine variety such that $S = \oplus_{i \geq 0}S_i$ ($S_0 = \mathbf{C}$) is a graded $\mathbf{C}$-algebra generated by $S_1$. Take the normalization $X$ of $V$ and let $R$ be its coordinate ring. Then the $\mathbf{C}^*$-action on $V$ extends to a $\mathbf{C}^*$-action on $X$. Assume that $X_{reg}$ admits a symplectic form $\omega$ and $(X, \omega)$ is a conical symplectic variety with $wt(\omega) = 1$ for this $\mathbf{C}^*$-action. Then $S_1 \subset R_1$,  
but $S_1 \ne R_1$ in general; hence $S \ne R(1)$. By the proof of Proposition 1 we see that $W := \mathrm{Spec}R(1)$ is a nilpotent orbit closure of a complex semisimple Lie algebra, but $V$ is not necessarily so. A similar situation appeared in [B-F].      
\vspace{0.2cm}

We are now going to characterize a conical symplectic variety obtained as a finite 
cover of a nilpotent orbit closure in a complex semisimple Lie algebra. \vspace{0.2cm}

{\bf Proposition 2}. {\em  Let $X = \mathrm{Spec} R$ be a conical symplectic variety of dimenion $2d$ with $wt(\omega) = l$. Denote by $R(l)$ the image of the natural homomorphism $\phi: \oplus Sym^i(R_l) \to R$
Assume that $R$ is a finite $R(l)$-module. 
Then } 

(i) {\em $R_l$ is a semisimple Lie algebra $\mathfrak{g}$ and the simply connected algebraic group $G$ with $\mathrm{Lie}(G) = \mathfrak{g}$ acts on $X$ as $\mathbf{C}^*$-equivariant symplectic automorphism.}  

 (ii) {\em The inclusion $R(l) \to R$ induces a finite $G$-morphism $f: X \to \bar{O}$ onto a nilpotent orbit closure $\bar{O}$ of $\mathfrak{g}$. Moreover, $f$ is \'{e}tale in codimension one and $f^{-1}(O) \to O$ is an \'{e}tale $G$-covering. The symplectic form $\omega$ is obtained as 
the pull-back of the Kirillov-Kostant form $\omega_{KK}$ on $O$ by $f$.} 
\vspace{0.2cm}

{\em Proof}. Let $\tilde{R(l)}$ be the integral closure of $R(l)$ in its quotient field and 
put $Y := \mathrm{Spec} \tilde{R(l)}$ and $Z: = \mathrm{Spec} R(l)$. Then the inclusions $R(l) \subset \tilde{R(l)} \subset R$ induces a 
finite surjective morphism $f: X \stackrel{\pi}\to Y \to Z$. The Poisson structure on $X$ 
naturally induces a Poisson structure on $Z$ because it has degree $-l$. By [Ka 2] it 
also induces a Poisson structure on the normalization $Y$ of $Z$. We prove that $\pi$ is 
\'{e}tale in codimension one. If not, then there is a smooth open subset $Y_0$ of $Y$ such that 
$\pi_0 (= \pi \vert_{\pi^{-1}(Y_0)}) : \pi^{-1}(Y_0) \to Y_0$ is branched over an irreducible smooth divisor $E$ on $Y_0$. Put $D := \pi^{-1}(E)$ and take a ramification point $p \in D$.   
One can take (analytic) local coordinates $(x, y_1, ... y_{2d-1})$  
of $X$ at $p$, and local coordinates $(z, y_1, ..., y_{2d-1})$ of $Y$ at $\pi (p)$ so that 
$\pi$ is defined by $z = x^m$ for the ramification index $m > 1$. Now let us consider the Poisson 2-vector $\theta$ on $X$. By the assumption $\theta$ descends to a Poisson 2-vector $\theta_Y$ on $Y$. Such a 2-vector that 
descends to a 2-vector on $Y$ can be locally written (around $p$) as follows.    
$$ \theta = \Sigma a_i(z, y_1, ..., y_{2d -1}) x \partial_x \wedge \partial_{y_i} 
+ \Sigma b_{i,j}(z, y_1, ..., y_{2d-1}) \partial_{y_i} \wedge \partial_{y_j}.$$ 
In particular, the $d$-th wedge product $\wedge^d \theta$ must be degenerate along $D$. 
This contradicts that $\wedge^d \theta$ is nowhere vanishing. Hence $\pi$ must be \'{e}tale 
in codimension one and $\theta_Y$ is non-degenerate 
on the regular locus $Y_{reg}$ of $Y$. The Poisson 2-vector $\theta_Y$ determines a symplectic 2-form $\omega_Y$ on $Y_{reg}$. 
In particular, $\wedge^d \omega_Y$ 
determines a nowhere vanishing section of $K_{Y_{reg}}$; hence it extends to a nowhere vanishing section of $K_Y$ and $K_Y = 0$. Since $\pi$ is \'{e}tale in codimension one, we have $K_X = \pi^*K_Y$. 
As $X$ is a symplectic variety, $X$ has only canonical singularities. We prove that $Y$ has 
canonical singularities. We take resolutions $\mu : \tilde{X} \to X$ and $\nu: \tilde{Y} \to Y$ 
so that the following diagram commutes: 

\begin{equation} 
\begin{CD} 
\tilde{X}  @>{\tilde{\pi}}>> \tilde{Y} \\ 
@V{\mu}VV @V{\nu}VV \\ 
X @>{\pi}>>  Y    
\end{CD} 
\end{equation} 

Let $\{E_i\}$ be the set of all $\nu$-exceptional prime divisors on $\tilde{Y}$. We write 
$K_{\tilde Y} = \nu^*K_Y + \Sigma a_iE_i$. Since $K_Y$ is Cartier, all $a_i$ are 
integers. We will get a contradiction by assuming that $a_{i_0} < 0$ for some $i_0$. 
Let $\{F_j\}$ be the set of all $\mu$-exceptional prime divisors on $\tilde{X}$.  
Write $\tilde{\pi}^*E_{i_0} = \Sigma r_j F_j$ for some integers $r_j$. (Note that every 
irreducible component of $\tilde{\pi}^{-1}(E_{i_0})$ is $\mu$-exceptional.)
There is an irreducible component of $\pi^{-1}(E_{i_0})$ which dominates $E_{i_0}$ 
by the map $\tilde{\pi}$. Assume that $F_{j_0}$ is such a component. Then  
$r_{j_0}$ is the ramification index of $\tilde{\pi}$ along $F_{j_0}$. 
Now we have $\tilde{\pi}^*K_{\tilde Y} = \tilde{\pi}^*\nu^*K_Y + \Sigma a_i \tilde{\pi}^*E_i$ 
On the other hand, we can write $K_{\tilde X} = \tilde{\pi}^*K_{\tilde Y} + \Sigma b_j F_j $ 
for some integers $b_j$.
Notice again that if $F_j$ is not $\tilde{\pi}$-exceptional, then $b_j = r_j - 1$ for the 
ramification index $r_j$ of $\tilde{\pi}$ along $F_j$. 
 By putting these together and using the fact $\tilde{\pi}^*\nu^*K_Y 
= \nu^*\pi^*K_Y = \mu^*K_X$, we finally get $$K_{\tilde X} 
=\mu^*K_Y + \Sigma a_i \tilde{\pi}^*E_i + \Sigma b_j F_j. $$  
Then the coefficients of $F_{j_0}$ on the right hand side equals 
$a_{i_0}r_{j_0} + r_{j_0} - 1$, which is negative becasue $a_{j_0}$ is a negative integer. 
This contradicts that $X$ has canonical singularities. Therefore all $a_i$ are non-negative   
and $Y$ has canonical singularities.   
This means that $(Y, \omega_Y)$ is a symplectic variety. The $\mathbf{C}^*$-action on $Z$ uniquely extend to the $\mathbf{C}^*$-action on $Y$, and $(Y, \omega_Y)$ becomes a conical symplectic variety.   
    
Let $\mathrm{Aut}^{{\mathbf C}^*}(X, \omega)$ be the group of $\mathbf{C}^*$-equivariant symplectic automorphisms of $(X, \omega)$. Then $\mathrm{Aut}^{{\mathbf C}^*}(X, \omega)$ acts on the graded 
ring $R$; hence acts on $R(l)$. Therefore $\mathrm{Aut}^{{\mathbf C}^*}(X, \omega)$ acts on $Z$. 
As $Y$ is the normalization of $Z$, it acts also on $Y$. 
Since $\mathrm{Aut}^{{\mathbf C}^*}(X, \omega)$ preserves $(\pi^{reg})^*\omega_Y$, the action 
of $\mathrm{Aut}^{{\mathbf C}^*}(X, \omega)$ on $Y$ preserves $\omega_Y$. This means that there is 
a homomorphism $\mathrm{Aut}^{{\mathbf C}^*}(X, \omega) \to \mathrm{Aut}^{{\mathbf C}^*}(Y, \omega_Y)$ of 
algebraic groups.  It induces a homomorphism $(\mathrm{Aut}^{{\mathbf C}^*}(X, \omega))^0 \to (\mathrm{Aut}^{{\mathbf C}^*}(Y, \omega_Y))^0$ of the identity components. 
The tangent space $T_{[id]}\mathrm{Aut}^{{\mathbf C}^*}(X, \omega)$ is isomorphic to  
$R_{\mathrm{wt}(\omega)} = R_l$ (cf. [Na-Part I, Proposition 2]). Similarly, $T_{[id]}\mathrm{Aut}^{{\mathbf C}^*}(Y, \omega_Y) 
\cong R(l)_{\mathrm{wt}(\omega_Y)} = R_l$. Therefore the homomorphism    
$(\mathrm{Aut}^{{\mathbf C}^*}(X, \omega))^0 \to (\mathrm{Aut}^{{\mathbf C}^*}(Y, \omega_Y))^0$ is a surjection 
with finite kernel. 
   
Now look at the normalization map $Y \to Z$.  
Since the degree of a homogeneous element of $R(l)$ is a multiple of $l$, the same is true for $\tilde{R(l)}$. In fact, the given $\mathbf{C}^*$-action on $Z$ factorizes as $\mathbf{C}^* \stackrel{\chi_l}\to \mathbf{C}^* \stackrel{\sigma}\to \mathrm{Aut}(Z)$. Here $\chi_l$ is defined by $\chi_l(t) = t^l$.  
Note that $R_l$ coincides with degree $1$ part of the new action $\sigma$. 
Since $Y$ is the normalization of $Z$, the action $\sigma$ extends to 
a $\mathbf{C}^*$-action, say $\tilde{\sigma}$, on $Y$. Then the action $\tilde{\sigma} \circ 
\chi_l : \mathbf{C}^* \to \mathrm{Aut}(Y)$ coincides with the original $\mathbf{C}^*$-action 
on $Y$. Thus we can write $\tilde{R(l)} = 
\oplus_{j \geq 0} (\tilde{R(l)})_{jl}$. We introduce a new grading of $\tilde{R(l)}$ by 
putting $(\tilde{R(l)})_{[j]} := \tilde{R(l)}_{jl}$. Then $R(l)$ is a subring of $\tilde{R(l)}$ generated 
by degree one elements. 

Now the situation is the same as in Proposition 1. Since $wt(\omega_Y) = 1$ for the new 
grading, we see that $\mathfrak{g} := R_l$ is a semisimple Lie algebra, $(\mathrm{Aut}^{{\mathbf C}^*}(Y, \omega_Y))^0$ is the adjoint group of $\mathfrak{g}$ and $Z$ is a nilpotent orbit closure  $\bar{O}$ of $\mathfrak{g}$. 

By the construction $f: X \to \bar{O}$ is a finite $(\mathrm{Aut}^{{\mathbf C}^*}(X, \omega))^0$- 
morphism. If $G$ is a simply connected semisimple algebraic group with 
$Lie(G) = \mathfrak{g}$, then there is a surjective homomorphism $G \to (\mathrm{Aut}^{{\mathbf C}^*}(X, \omega))^0$ with finite kernel. Hence $f$ is a finite $G$-morphism. Since $\mathrm{Codim}( \mathrm{Sing}(\bar{O})) \geq 2$, the normalization map is \'{e}tale in codimension one. We have 
already seen that $\pi$ is \'{e}tale in codimension one; hence, $f$ is \'{e}tale in codimension one.  
Take a point $x \in f^{-1}(O)$. Since $O$ is the $G$-orbit containing $f(x)$, the point $x$ is contained in a $G$-orbit $O'$ on $X$ such that $f(O') = O$. This implies that $O'$ is dense 
in $f^{-1}(O)$. But it holds for any other $y \in f^{-1}(O)$. This is possible only when $y \in O'$.  Hence we have $O' = f^{-1}(O)$.    
Q.E.D. \vspace{0.2cm}

Let $G$ be a simply connected complex semisimple algebraic group, and let $(X, \omega)$ be an affine  symplectic variety obtained as a finite $G$-covering $f: X \to \bar{O}$ of a nilpotent 
orbit closure $\bar{O}$ in $\mathfrak{g}$ where $f^{-1}(O) \to O$ is an \'{e}tale G-covering. Moreover we assume that $\omega$ is the pull-back of the pull-back of the Kirillov-Kostant form $\omega_{KK}$ on $O$ by $f$ Such a variety $X$ has been studied by Brylinski and Kostant [B-K]. 
The nilpotent orbit $O$ has a $\mathbf{C}^*$-action 
induced by the scalar $\mathbf{C}^*$-action on $\mathfrak{g}$. This $\mathbf{C}^*$-action 
is generated by the Euler vector field $\bar{\eta} = \Sigma x_i \partial_{x_i}$ where $x_i$ are 
basis of ${\mathfrak{g}}^*$. \vspace{0.15cm}

{\bf Lemma} ([B-K, \S 1]). {\em ${f}^{-1}(O)$ has a $\mathbf{C}^*$-action generated by a vector field $\eta$ with ${f}_*\eta = 2\bar{\eta}$.}

{\em Proof}. Take $e \in O$ and choose $\epsilon \in {f}^{-1}(O)$ so that $f(\epsilon) = 
e$. Let $G^e$ and $G^{\epsilon}$ be the stabilizer groups of $e$ and $\epsilon$. 
Then $({G}^e)^0 \subset {G}^{\epsilon} \subset {G}^e$, where $({G}^e)^0$ is the identity 
component of ${G}^e$. In particular, we have ${\mathfrak{g}}^e = {\mathfrak{g}}^{\epsilon}$.  
Let $N^{\epsilon}$ be the normalizer of ${G}^{\epsilon}$ in $G$. 
The group $N^{\epsilon}/{G}^{\epsilon}$ acts on ${f}^{-1}(O) = G/{G}^{\epsilon}$ by 
$\bar{n}\cdot g{G}^{\epsilon} := gn^{-1}{G}^{\epsilon}$ for $n \in N^{\epsilon}$ and $g \in G$. 
By definition this action commutes with the $G$-action on ${f}^{-1}(O)$. 

By the Jacobson-Morozov theorem we can take elements $f$ and $h$ of $\mathfrak{g}$ 
so that $[h, e] = 2e$, $[h, f] = -2f$ and $[e, f] = h$. Then we can show that $h \in \mathfrak{n}^{\epsilon}$. 
In fact, for $x \in {\mathfrak{g}}^e$, we have 
$$ [[h, x], e] = -[[x, e], h] - [[e,h], x] = 0, $$ which implies that $[h, x] \in {\mathfrak{g}}^e$. 
Since ${\mathfrak{g}}^e = {\mathfrak{g}}^{\epsilon}$, we see that 
$[h, {\mathfrak{g}}^{\epsilon}] \subset {\mathfrak{g}}^{\epsilon}$. 
Therefore, $\mathrm{exp}(\mathbf{C}h) \subset 
N^{\epsilon}$ so that $\mathbf{C}^* = \mathrm{exp}(\mathbf{C}h)$ defines a subgroup $C$ of $N^{\epsilon}/{G}^{\epsilon}$. The surjective homomorphism $\mathbf{C}^* \to C$ determines 
a $\mathbf{C}^*$-action on $X$. The final assertion follows from $[h, e] = 2e$.  Q.E.D
\vspace{0.15cm}

The $\mathbf{C}^*$-action uniquely extends to a $\mathbf{C}^*$-action on $X$ because $X$ is a normal variety with $\mathrm{Codim}_X(X - {f'}^{-1}(O)) \geq 2$. We call this $\mathbf{C}^*$-action the {\em right} $\mathbf{C}^*$-action. The right $\mathbf{C}^*$-action makes the coordinate ring $R$ a graded $\mathbf{C}$-algebra $R = \oplus R_{[i]}$  
The Poisson structure on $R$ has degree $-2$ with respect to this grading by the lemma above.  Let us denote by $R[2]$ the image of the natural homomorphism $\phi: \oplus Sym^i(R_{[2]}) \to R$; then the coordinate ring $\mathbf{C}[\bar{O}]$ of $\bar{O}$ is contained in $R[2]$. The coordinate ring $R$ is a finite $R[2]$-module.

Let us consider the conical symplectic variety $(X, \omega)$ and a finite surjective morphism $f: X \to \bar{O}$ in Proposition 2. The observations above can be applied to this $(X, \omega)$. Then $X$ has two $\mathbf{C}^*$-actions, one of them is the 
original one and another one is the {\em right} $\mathbf{C}^*$-action. Correspondingly $R$ 
has two gradings $R = \oplus R_i$ and $R = \oplus R_{[i]}$. 
\vspace{0.2cm}

%By the construction $R_l = R_{[2]}$ 
%and $\mathbf{C}[O] = R(l) = R[2]$. Assume that $l$ is an even (positive) integer, i.e.  
%$l = 2k$; then $R_{ik} = R_{[i]}$ for all $i$. Assume that $l$ is an odd integer; then 
%$R = \oplus R_{il}$, $R = \oplus R_{[2i]}$ and $R_{il} = R_{[2i]}$. 
%\vspace{0.2cm}

%{\bf Definition}. Let $R = \oplus_{i \geq 0} R_i$ be a graded ring. Assume that there is a positive integer $k$ such that $R_i = 0$ except when $i$ is some multiple of $k$. Then one can write 
%$R = \oplus_{j \ge 0}R_{jk}$. In such a case one can put $R_{[j]} := R_{jk}$ and introduce a new 
%grading $R = \oplus_{j \geq 0}R_{[j]}$. Such a grading is called a {\em reduced grading} of $R$.   
%\vspace{0.2cm}

{\bf Proposition 3}. {\em  Let $X = \mathrm{Spec} R$ be a conical symplectic variety with $wt(\omega) = l$. Denote by $R(l)$ the image of the natural homomorphism $\phi: \oplus Sym^i(R_l) \to R$. Assume that $R$ is a finite $R(l)$-module.  
When $l $ is an even number, one has $R = \oplus_{j \geq 0}R_{jl/2}$.  When 
$l $ is an odd number, one has $R = \oplus_{j \geq 0}R_{jl}$.}
\vspace{0.15cm}

{\em Proof}.  The original $\mathbf{C}^*$-action and the right $\mathbf{C}^*$-action respectively determine homomorphism $\sigma: \mathbf{C}^* \to \mathrm{Aut}(X)$ and $\tau: \mathbf{C}^* \to \mathrm{Aut}(X)$. 

Assume that $l$ is even, i.e. $l = 2k$. Consider the homomorphism  $\chi_k : 
\mathbf{C}^* \to \mathbf{C}^*$ defined by $\chi_k(t) = t^k$. Then $\sigma$ and $\tau^k :=  \tau \circ 
\chi_k$ induce the same $\mathbf{C}^*$-action on $\bar{O}$ by the construction. 
This means $\sigma = \tau^k$. The degree $i$ part of $R$ with 
respect to $\sigma$ is $R_i$. On the other hand, the degree $i$ part of $R$ with respect to 
$\tau^k$ is, by definition, $R_{[i/k]}$. This means that $R_i \ne 0$ only if $i$ is some  
multiple of $k$. Namely, $R = \oplus_{i \geq 0}R_{il/2}$. 

Assume that $l$ is odd. Then we can take integers $p$ and $q$ so that $pl + 2q = 1$. 
Let us consider the $\mathbf{C}^*$-action $\sigma^p \circ \tau^q$ on $X$. 
By the choice of $p$ and $q$, it induces the scalar $\mathbf{C}^*$-action on 
$\bar{O}$. Let $R = \oplus_{i \geq 0} R_{<i>}$ be the grading determined by $\sigma^p \circ \tau^q$. Then $R_{<i>} = R_{il}$. Hence $R = \oplus_{i \geq 0}R_{il}$. Q.E.D.

Department of Mathematics, Graduate school of Science, Kyoto University 

e-mail: namikawa@math.kyoto-u.ac.jp 

\end{document}